\documentclass[11]{article}

\usepackage{enumerate}
\usepackage{amsmath}
\usepackage{amsfonts}
\usepackage{amsthm}

\def\ind{{\rm 1\kern-.27em I}}

\newcommand{\dis}{\displaystyle}
\newcommand{\eps}{\varepsilon}
\newcommand{\ud}{\mathrm{d}}
\newcommand{\D}{\mathbb{D}}
\newcommand{\C}{\mathbb{C}}

\newcommand{\N}{\mathbb{N}}
\newcommand{\T}{\mathbb{T}}

\newcommand{\CQFD}
{%
\mbox{}%
\nolinebreak%
\hfill%
\rule{2mm}{2mm}%
\medbreak%
\par%
}
\newtheorem{Theo}{Theorem}[section]
\newtheorem{Lem}[Theo]{Lemma}
\newtheorem{Prop}[Theo]{Proposition}
\newtheorem{Def}[Theo]{Definition}
\newtheorem{Cor}[Theo]{Corollary}

\newtheoremstyle{Remark}{3pt}{3pt}{}{}{}{.}{6pt}{}
\theoremstyle{Remark}
\newtheorem{Rem}[Theo]{\bf Remark}

\numberwithin{equation}{section}

\begin{document}
\sloppy
\title{\sc Weighted composition operators as Daugavet centers}
\author{R. Demazeux}
\date{}
\maketitle

\thispagestyle{empty}

\vspace{-25pt}

\begin{abstract}
We investigate the norm identity $\|uC_\varphi+T\|=\|u\|_\infty+\|T\|$ for classes of operators on $C(S)$, where $S$ is a compact Hausdorff space without isolated point, and characterize those weighted composition operators which satisfy this equation for every weakly compact operator~$T:C(S)\rightarrow C(S).$ We also give a characterization of such weighted composition operator acting on the disk algebra~$A(\D).$
\end{abstract}

\smallskip

\setcounter{Theo}{0}\setcounter{section}{0}

\section{Introduction}
In $1963,$ Daugavet proved~\cite{D} the norm equality
\begin{equation}\label{eqdaug}
\|I+T\|=1+\|T\|
\end{equation}
now known as the Daugavet equation for every compact operator $T:C([0,1])\rightarrow C([0,1]).$ Over the years, this property was extended to larger classes of operators and various spaces~: $C(S)$ where $S$ is a compact Hausdorff space without isolated point~\cite{FS}, $L^1(\mu)$ for measure~$\mu$ without atoms~\cite{Lo}, the disk algebra~$A(\D)$ or the Hardy space~$H^\infty$~\cite{Wo}. Actually, if~\eqref{eqdaug} holds for every rank-$1$ operators on~$X$, then it holds for every weakly compact operators, $X$ contains a copy of~$\ell_1$, $X$ cannot have an unconditional basis and even cannot embed into a space having an unconditional basis (see the survey~\cite{We2}).\par
Recently in~\cite{P}, the author showed that if we substitute the Identity in~\eqref{eqdaug} with an into isometry~$J:L^1[0,1]\rightarrow L^1[0,1]$ then equation $\|J+T\|=1+\|T\|$ holds for narrow operators, and particularly for weakly compact operators on~$L^1[0,1].$ In arbitrary Banach spaces, this has been investigated by T. Bosenko and V. Kadets in~\cite{BK}. They introduced the following concept~:
\begin{Def}
Let $X$ be a Banach space. A linear continuous operator $G:X\rightarrow X$ is said to be a Daugavet center if the norm identity 
\begin{equation}\label{daugcenter}
\|G+T\|=\|G\|+\|T\|
\end{equation}
holds for every rank-$1$ operator~$T:X\rightarrow X.$
\end{Def}
Bosenko and Kadets showed that if $G:X\rightarrow X$ is a non zero Daugavet center then equation~\eqref{daugcenter} holds true for every strong Radon-Nikod\'ym operator on~$X$, and so for weakly compact operators on~$X.$
Moreover $G$ fixes a copy of~$\ell_1,$ and $X$ cannot have an unconditional basis, merely $X$ cannot be embedded into a space having an unconditional basis.\par
In section~$2$ of this paper, we give a characterization of weighted composition operators on~$C(S)$ which are Daugavet centers ($S$ compact Hausdorff space without isolated point). We give examples of Daugavet centers which are not weighted composition operators and prove that the set of Daugavet centers in~$C(S)$ is not convex. We also study equation~\eqref{daugcenter} for the class of operators whose adjoint has separable range. This encompass the class of operators factorizing through~$c_0.$\par
In section~$3$, we adapt D. Werner's method showing that certain function spaces have the Daugavet property (meaning that the identity operator is a Daugavet center) to characterize weighted composition operators on the disk algebra~$A(\D)$ which are Daugavet centers.

\section{Weighted composition operators as Daugavet center in $C(S)$}
Let $S$ denotes a compact Hausdorff space without isolated point. Considering continuous maps $\varphi:S\rightarrow S$ and $u:S\rightarrow \C$, we study the weighted composition operator $uC_\varphi:C(S)\rightarrow C(S)$ defined by $uC_\varphi(f)=u.(f\circ\varphi)$ for all $f\in C(S).$ We clearly have $\|uC_\varphi\|=\|u\|_\infty.$
We investigate the following equation~:\[\|uC_\varphi+T\|=\|u\|_\infty+\|T\| \qquad\qquad (E_{u,\varphi})\]
Note that if we take $\varphi(s)=s$ (for all $s$ in $S)$ and $u$ the constant function equal to~$1,$ the previous equation becomes the classical Daugavet equation. We will suppose that $u$ is not the constant function equal to zero. We want to find conditions on $\varphi$ and $u$ implying that every weakly compact operators on C(S) satisfies equation~$(E_{u,\varphi}).$ A first remark is that $u$ and $\varphi$ must be such that the operator $uC_\varphi$ is not itself compact. By a result of Kamowitz~\cite{Ka}, $uC_\varphi$ is compact if and only if $\varphi$ is constant on a neighborhood of each connected component of the set where $u$ is nonzero. So $\varphi$ should be non constant over at least one nonempty open set in $S.$\par
The main result of this section is the following theorem~:
\begin{Theo}
Let $S$ be a compact Hausdorff space without isolated point, $u\in C(S)$ and $\varphi$ be a continuous function from~$S$ to~$S.$\\Then every weakly compact operators $T:C(S)\rightarrow C(S)$ satisfies equation $(E_{u,\varphi})$ if and only if $\varphi^{-1}(\{t\})$ is nowhere dense in~$S,$ for every $t\in S$ and $\vert u\vert$ is constant on $S.$
\end{Theo}
The straightforward direction was already proved in~\cite{BK} for rank-$1$ operators and for~\mbox{$u\equiv 1.$} Here we give a direct and simple proof for weakly compact operators, and we check that conditions on~$\varphi$ and~$u$ are necessary.\par
We first begin with some notations and terminology. The dual space of $C(S)$ consisting of all regular borel measures on $S$ of finite variation will be denoted by $M(S).$ If $s\in S,$ we define the corresponding Dirac functional $\delta_s$ by $\delta_s(f)=f(s)$ for every $f\in C(S).$ Then $\delta_s\in M(S)$ and $\|\delta_s\|=1.$\par

Following~\cite{We1}, the key idea is to represent an operator $T:C(S)\to C(S)$ by the family of measures $(\mu_s)_{s\in S}$ defined by $\mu_s=T^*(\delta_s)$, such that :\[(Tf)(s)=\left\langle Tf,\delta_s\right\rangle=\left\langle f,\mu_s\right\rangle=\int_Sf\ \mathrm{d}\mu_s.\]
Thus, the (weakly) compact nature  of $T$ is reformulated in terms of continuity of the map $s\mapsto \mu_s$ in the following (c.f.~\cite{DS}, Th. VI, 7.1)~:  
\begin{Lem}\label{cont}
Let $T:C(S)\to C(S)$ be an operator and $(\mu_s)_{s\in S}$ be the family of measures associated to~$T.$ Then : 
\begin{enumerate}[i)]
  \item $s\mapsto\mu_s$ is continuous from $S$ to $M(S)=C(S)^*$ endowed with the weak$^*$-topology, i.e. $\sigma(M(S),C(S))$.
  \item $T$ is weakly compact if and only if $s\mapsto\mu_s$ is continuous for the weak-topology on $M(S)$, \textit{i.e.} for $\sigma\big(M(S),M(S)^*\big)$.
  \item $T$ is compact if and only if $s \mapsto \mu_s$ is continuous for the norm topology on $M(S)$.
\end{enumerate}
\end{Lem}

Note that $\|T\|=\dis\sup_{s\in S} \|\mu_s\|$, and that the operator $uC_\varphi$ is represented by the family of measures $(u(s)\delta_{\varphi(s)})_{s\in S}.$ Indeed~: \[(uC_\varphi)^*(\delta_s)(f)=\delta_s(u.f\circ\varphi)=u(s)f\big(\varphi(s)\big)=u(s)\delta_{\varphi(s)}(f).\]

The following proposition shows, assuming $\vert u\vert$ is constant, that for every operator $T$ on $C(S),$ there is a $\lambda \in\T$ such that $\lambda T$ satisfies equation $(E_{u,\varphi})$.
\begin{Prop}
Let $S$ be a compact Hausdorff space, and $T:C(S)\rightarrow C(S)$ be an operator. Assume that $\vert u\vert$ is constant. Then \[\max_{\lambda \in\T}{\|uC_\varphi+\lambda T\|}=\|u\|_\infty+\|T\|.\]
\end{Prop}
{\bf Proof.} Let $(\mu_s)_{s\in S}$ be the family of measures associated to $T$. Then
\begin{align*}
\max_{\lambda \in\T} \|uC_\varphi+\lambda T\| & = \max_{\lambda \in\T} \sup_{s\in S}\|u(s)\delta_{\varphi(s)}+\lambda \mu_s\|\\
                    & = \sup_{s\in S}\max_{\lambda \in\T} \Big(\vert u(s)\delta_{\varphi(s)}+\lambda \mu_s\vert\big(\{\varphi(s)\}\big)+\vert u(s)\delta_{\varphi(s)}+\lambda \mu_s\vert\big(S\backslash\{\varphi(s)\}\big)\Big) \\
                    & = \sup_{s\in S}\max_{\lambda \in\T} \Big(\big\vert u(s)+\lambda \mu_s\big(\{\varphi(s)\}\big)\big\vert+\vert\mu_s\vert\big(S\backslash\{\varphi(s)\}\big)\Big) \\ 
                    & = \sup_{s\in S}\Big(\vert u(s)\vert+\big\vert\mu_s\big(\{\varphi(s)\}\big)\big\vert+\vert\mu_s\vert\big(S\backslash\{\varphi(s)\}\big)\Big) \\
                    & = \sup_{s\in S}(\|u\|_\infty+\|\mu_s\|)\qquad\quad\textrm{since }\vert u(s)\vert=\|u\|_\infty\\
                    & = \|u\|_\infty+\|T\|.             
\end{align*}
\CQFD
\begin{Rem}
  \begin{enumerate}[i)]
     \item In the real case, a similar result holds replacing $``\lambda\in\T$'' by $``\lambda\in\{\pm 1\}$''.
     \item Without assumption on the modulus of $u$, the previous result is not true anymore. For instance, taking $v\in C(S),$ we have that $\max_{\lambda\in\T}\|uC_\varphi+\lambda vC_\psi\|=\|\vert u\vert+\vert v\vert\|_\infty$ which is not equal to $\|u\|_\infty+\|v\|_\infty$ in general.
   \end{enumerate}
\end{Rem}
\subsection{Equation $(E_{u,\varphi})$ for weakly compact operators on $C(S)$}
Let $T:C(S)\rightarrow C(S)$ be an operator and $(\mu_s)_{s\in S}$ be the family of measures associated to~$T.$ Then 
\[\|uC_\varphi+T\|=\sup_{s\in S}\|u(s)\delta_{\varphi(s)}+ \mu_s\|=\sup_{s\in S}\Big(\big\vert u(s)+ \mu_s\big(\{\varphi(s)\}\big)\big\vert+\vert\mu_s\vert\big(S\backslash\{\varphi(s)\}\big)\Big)\]
and \[\|u\|_\infty+\|T\|=\sup_{s\in S}(\|u\|_\infty+\|\mu_s\|)=\sup_{s\in S}\Big(\|u\|_\infty+\big\vert\mu_s\big(\{\varphi(s)\}\big)\big\vert+\vert\mu_s\vert\big(S\backslash\{\varphi(s)\}\big)\Big).\]

We have the following lemma which gives a characterization of the operators satisfying equation $(E_{u,\varphi})~:$ 

\begin{Lem}\label{eq1}
\[\|uC_\varphi+T\|=\|u\|_\infty+\|T\|\] if and only if  
\begin{equation}
\dis\sup_{\{s\in S\mid\ \|\mu_s\|>\|T\|-\eps\}}\Big(\big\vert u(s)+\mu_s\big(\{\varphi(s)\}\big)\big\vert-\big(\|u\|_\infty+\big\vert\mu_s(\{\varphi(s)\})\big\vert\big)\Big)=0 
\end{equation}
 for all $\eps>0.$
\end{Lem}
{\bf Proof.} Sufficient condition : 
let $\eps>0$ and $U=\{s\in S \mid \|\mu_s\|>\|T\|-\eps\}$ which is not empty. Then~: 
 \begin{align*}
 \|uC_\varphi+T\| & \geq \sup_{s\in U}\|u(s)\delta_{\varphi(s)}+\mu_s\| \\
            & \geq \sup_{s\in U}\Big(\big\vert u(s)+\mu_s\big(\{\varphi(s)\}\big)\big\vert + \vert\mu_s\vert\big(S\backslash\{\varphi(s)\}\big)\Big) \\
            & \geq \sup_{s\in U}\bigg(\big\vert u(s)+\mu_s\big(\{\varphi(s)\}\big)\big\vert+\|u\|_\infty+\|\mu_s\|-\Big(\|u\|_\infty+\big\vert\mu_s\big(\{\varphi(s)\}\big)\big\vert\Big)\bigg)\\
            & \geq \|u\|_\infty+\|T\|-\eps+\sup_{s\in U}\bigg(\big\vert u(s)+\mu_s\big(\{\varphi(s)\}\big)\big\vert-\Big(\|u\|_\infty+\big\vert\mu_s\big(\{\varphi(s)\}\big)\big\vert\Big)\bigg) \\
            & \geq \|u\|_\infty+\|T\|-\eps \qquad \big(\textrm{with}\ (2.1)\big).
 \end{align*}\\
Necessary condition : let us assume that there exist $\alpha$ and $\eps>0$ such that for all $s\in S$, $\|\mu_s\|>\|T\|-\eps$ implies \[\big\vert u(s)+\mu_s\big(\{\varphi(s)\}\big)\big\vert-\Big(\|u\|_\infty+\big\vert\mu_s\big(\{\varphi(s)\}\big)\big\vert\Big)<-\alpha<0.\] Then
 \begin{align*}
 \|uC_\varphi+T\| & = \sup_{s\in S}\|u(s)\delta_{\varphi(s)}+\mu_s\| \\
            & = \max{\Big(\sup_{\{s \mid\ \|\mu_s\|>\|T\|-\eps\}}\|u(s)\delta_{\varphi(s)}+\mu_s\|,\sup_{\{s \mid\ \|\mu_s\|\leq\|T\|-\eps\}}\|u(s)\delta_{\varphi(s)}+\mu_s\|\Big)}. \\ 
 \end{align*}
The second term is lower than $\|u\|_\infty+\|T\|-\eps$. For the first term, we write as before
\begin{align*}
 &\sup_{\{s\in S\mid\ \|\mu_s\|>\|T\|-\eps\}}\|u(s)\delta_{\varphi(s)}+\mu_s\| \\
 &= \sup_{\{s\in S\mid\ \|\mu_s\|>\|T\|-\eps\}}\Big(\big\vert u(s)+\mu_s\big(\{\varphi(s)\}\big)\big\vert+\vert\mu_s\vert\big(S\backslash\{\varphi(s)\}\big)\Big)\\
 &          = \sup_{\{s\in S\mid\ \|\mu_s\|>\|T\|-\eps\}}\bigg(\big\vert u(s)+\mu_s\big(\{\varphi(s)\}\big)\big\vert+\|u\|_\infty+\|\mu_s\|-\Big(\|u\|_\infty+\big\vert\mu_s\big(\{\varphi(s)\}\big)\big\vert\Big)\bigg) \\
 &            \leq \|u\|_\infty+\|T\|+\sup_{\{s\in S\mid\ \|\mu_s\|>\|T\|-\eps\}}\Big(\big\vert u(s)+\mu_s\big(\{\varphi(s)\}\big)\big\vert-\big(\|u\|_\infty+\big\vert\mu_s(\{\varphi(s)\})\big\vert\big)\Big)\\
 &            \leq \|u\|_\infty+\|T\|-\alpha.
\end{align*}
Thus $\|uC_\varphi+T\|<\|u\|_\infty+\|T\|-\min{(\eps,\alpha)}<\|u\|_\infty+\|T\|$, which leads to a contradiction.
\CQFD
As a consequence, we state the following useful corollary~:
\begin{Cor}\label{eq2}
Assume that the family $(\mu_s)_{s\in S}$ satisfies the following condition~: 
for every nonempty open set$\ U\subset S,$
\begin{equation}
\sup_{s\in U}\bigg(\big\vert u(s)+\mu_s\big(\{\varphi(s)\}\big)\big\vert-\Big(\|u\|_\infty+\big\vert\mu_s\big(\{\varphi(s)\}\big)\big\vert\Big)\bigg)=0.
\end{equation}
Then \[\|uC_\varphi+T\|=\|u\|_\infty+\|T\|.\]
\end{Cor}
{\bf Proof.} Take $\eps>0$, and call $U=\{s\in S \mid\ \|\mu_s\|>\|T\|-\eps\}.$ Thanks to Lemma~\ref{eq1}, we only need to show that $U$ is a nonempty open subset of~$S$. It is clear that~$U$ is nonempty. Take $s_0\in U.$ There exists $f_0\in C(S),\ \|f_0\|_\infty\leq1$ such that $\vert \mu_{s_0}(f_0)\vert>\|T\|-\eps.$ From Lemma~\ref{cont}, we know that $s\mapsto\mu_s$ is continuous for the weak$^*$-topology on $M(S),$ hence $s\mapsto\mu_s(f_0)$ is continuous. Then $V=\{s\in S\mid\ \vert\mu_s(f_0)\vert>\|T\|-\eps\}$ is an open neighborhood of $s_0$ contained in $U.$ So $U$ is a nonempty open subset of~$S.$
\CQFD

Now we can show a first result dealing with weakly compact operators. The following theorem gives sufficient conditions on $u$ and $\varphi$ implying that every weakly compact operator on $C(S)$ satisfies equation $(E_{u,\varphi}).$
 
\begin{Theo}\label{faiblcomp}
Let $S$ be a compact Hausdorff space (without isolated point). Assume that $\vert u\vert$ is constant on~$S$ and $\varphi(U)$ is infinite for every nonempty open subset $U$ of $S.$\\ 
Then $\|uC_\varphi+T\|=\|u\|_\infty+\|T\|$ for every weakly compact operator $T:C(S)\rightarrow C(S).$ 
\end{Theo}
Note that condition on $\varphi$ forces $S$ to have no isolated point.\\
{\bf Proof.} Assume that $\varphi(U)$ is infinite for every nonempty open subset $U$ of $S$ and that $\vert u\vert$ is constant. If the family of measures $(\mu_s)_{s\in S}$ representing $T$ does not satisfy~$(2.2)$ of Corollary~\ref{eq2}, then there exist a nonempty open set~$U\subset S$ and $\beta>0$ such that \[\big\vert u(s)+\mu_s\big(\{\varphi(s)\}\big)\big\vert-\Big(\|u\|_\infty+\big\vert\mu_s\big(\{\varphi(s)\}\big)\big\vert\Big)<-2\beta\,\quad \forall s\in U.\]In particular we have, since $\vert u(s)\vert=\|u\|_\infty$ for all $s\in S :$
\begin{align*}
\big\vert\mu_s\big(\{\varphi(s)\}\big)\big\vert & >2\beta-\|u\|_\infty+\big\vert u(s)+\mu_s\big(\{\varphi(s)\}\big)\big\vert\\
 & \geq2\beta-\|u\|_\infty+\vert u(s)\vert-\big\vert\mu_s\big(\{\varphi(s)\}\big)\big\vert\\
 & =2\beta-\big\vert\mu_s\big(\{\varphi(s)\}\big)\big\vert
\end{align*} which gives \[\big\vert\mu_s\big(\{\varphi(s)\}\big)\big\vert>\beta,\quad\textrm{for all }s\in U.\]
Take $t\in S$. Then $s\in S\mapsto\mu_s(\{t\})\in\C $ is continuous. Indeed : from Lemma~\ref{cont}, $s\mapsto\mu_s$ is continuous for the weak-topology on $M(S)$. Since $\mu\mapsto\mu(\{t\})$ belongs to $M(S)^*$, it is continuous on $M(S)$ endowed with the weak-topology.\\
Let $s_0\in U$, and define \[U_1=\big\{s\in U\mid\ \big\vert\mu_s\big(\{\varphi(s_0)\}\big)\big\vert>\beta\big\}.\] From above, $U_1$ is an open subset of $U$ (and so of $S$) which contains $s_0$. Since $\varphi(U_1)$ is infinite, one can find $s_1$ in $U_1$ satifying $\varphi(s_1)\neq\varphi(s_0).$ Then we have
\begin{align*}
\big\vert\mu_{s_1}\big(\{\varphi(s_1)\}\big)\big\vert>\beta ,&\quad \textrm{since}\ s_1\in U\\ 
\big\vert\mu_{s_1}\big(\{\varphi(s_0)\}\big)\big\vert>\beta.&
\end{align*}
Consider now  \[U_2=\big\{s\in U_1\mid\ \big\vert\mu_s\big(\{\varphi(s_1)\}\big)\big\vert>\beta\big\}.\] It is an open subset of $U$ containing $s_1$, and it contains an element $s_2$ such that $\varphi(s_2)\neq\varphi(s_0)$ and $\varphi(s_2)\neq\varphi(s_1)$ (since $\varphi(U_2)$ is infinite). Then we have, since $s_2\in U_2\subset U_1\subset U,$
\begin{align*}
\big\vert\mu_{s_2}\big(\{\varphi(s_2)\}\big)\big\vert>\beta&\\
\big\vert\mu_{s_2}\big(\{\varphi(s_1)\}\big)\big\vert>\beta&\\ 
\big\vert\mu_{s_2}\big(\{\varphi(s_0)\}\big)\big\vert>\beta&.
\end{align*}
In such a way we construct a decreasing sequence of open subsets $U_n\subset U$, and a sequence of elements $(s_n)_{n\geq0}$, $s_n\in U_n$ having the property
\begin{align*}
U_{n+1} & =  \big\{s\in U_n\mid\ \big\vert\mu_s\big(\{\varphi(s_n)\}\big)\big\vert>\beta\big\}, \\
s_{n+1} & \in  U_{n+1}\\
\varphi(s_{n+1}) & \notin\{\varphi(s_0),\ldots,\varphi(s_n)\}.
\end{align*}
So 
\[\big\vert\mu_{s_n}\big(\{\varphi(s_j)\}\big)\big\vert>\beta,\qquad j=0,\ldots,n-1\]
which leads to a contradiction writing that \[\|T\|\geq\|\mu_{s_n}\|\geq\vert\mu_{s_n}\vert\big(\{\varphi(s_0),\ldots,\varphi(s_{n-1})\}\big)\geq n\beta,\qquad \forall n\in\N.\]
\CQFD
We now give necessary conditions on $\varphi$ and $u$ to ensure that every weakly compact operator on $C(S)$ satisfies equation $(E_{u,\varphi}).$ Actually we only need to consider rank-$1$ operators. 
\begin{Theo}\label{rec1}
Let $S$ be a compact Hausdorff space without isolated point. Assume that every rank-$1$ operator on $C(S)$ satisfies equation $(E_{u,\varphi}).$ Then $\vert u\vert$ is constant and $\varphi^{-1}(\{t\})$ is nowhere dense in $S$, for every $t\in S.$
\end{Theo}
{\bf Proof : }We first show that $\vert u\vert$ is constant on~$S.$ Arguing by contradiction, assume there exists $s_0\in S$ such that $\vert u(s_0)\vert<\|u\|_\infty.$ Then there exists $\delta>0$ and an open neighborhood $U$ of $s_0$ satisfying \[\forall s\in U,\ \vert u(s)\vert<\|u\|_\infty-\delta.\]
Choose a continuous function $v$ such that~:~$0\leq v\leq1,\ v(s_0)=1$ and $v(s)<1$ for all $s\neq s_0.$ We define the operator $T=v\delta_\tau$ where $\tau$ is an element of ~$S.$ Then $\mu_s=T^*(\delta_s)=v(s)\delta_\tau,\ \|\mu_s\|=v(s)$. Choose $\eps>0$ such that we have $\{s\in S\mid\ v(s)>1-\eps\}\subset U.$ It follows that
\begin{align*}
&\sup_{\{s \mid\ v(s)>1-\eps\}}\bigg(\big\vert u(s)+v(s)\delta_\tau\big(\{\varphi(s)\}\big)\big\vert-\Big(\|u\|_\infty+\vert v(s)\vert\delta_\tau\big(\{\varphi(s)\}\big)\Big)\bigg)\\
&\leq\sup_{s\in U}\bigg(\vert u(s)\vert+v(s)\delta_\tau\big(\{\varphi(s)\}\big)-\Big(\|u\|_\infty+v(s)\delta_\tau\big(\{\varphi(s)\}\big)\Big)\bigg) \\
&\leq\sup_{s\in U} \vert u(s)\vert-\|u\|_\infty\\
&\leq-\delta<0.
\end{align*}
The family of measures $(\mu_s)_{s\in S}$ does not satisfy condition $(2.1)$ of Lemma~\ref{eq1} and consequently $T$ does not satisfy equation $(E_{u,\varphi}),$ which is false since $T$ is a rank-$1$ operator. So $\vert u\vert$ is constant.\par

Now we prove that for every $t\in S,\ \varphi^{-1}(\{t\})$ is nowhere dense in $S.$ Let $U$ be a nonempty open subset of~$S.$ We want to find $s\in U$ such that $\varphi(s)\neq t.$ Consider the rank-$1$ operator $T=\delta_t gu,$ where $g\in C(S)$ such that $-1\leq g\leq-\frac{1}{2},\ g=-\frac{1}{2}$ outside $U$ and $\|g\|_\infty=1.$ Then $\|T\|=\|u\|_\infty>0$, the family of measures associated to~$T$ is given by $\mu_s=T^*(\delta_s)=u(s)g(s)\delta_t,$ and $\|\mu_s\|=\vert u(s)g(s)\vert=\|u\|_\infty\vert g(s)\vert$. \\
Take $\eps=\|u\|_\infty/2$ so that $V=\{s\in S\mid\ \vert u(s)g(s)\vert>\frac{\|u\|_\infty}{2}\}\subset U.$ Since $T$ satisfies equation $(E_{u,\varphi}),$ the family of measures $(\mu_s)_{s\in S}$ satisfies condition $(2.1)$ of Lemma~\ref{eq1}~: 
\begin{align*}
0 &=\sup_V\bigg(\big\vert u(s)+u(s)g(s)\delta_t\big(\{\varphi(s)\}\big)\big\vert-\Big(\|u\|_\infty+\vert u(s)g(s)\vert\delta_t\big(\{\varphi(s)\}\big)\Big)\bigg)\\
  &=\sup_V\Big(2g(s)\|u\|_\infty\delta_t\big(\{\varphi(s)\}\big)\Big)
\end{align*}
which is less than $\dis\sup_{s\in U}\Big(2g(s)\|u\|_\infty\delta_t\big(\{\varphi(s)\}\big)\Big).$ It follows that there exists $s\in U$ such that $\delta_t\big(\{\varphi(s)\}\big)=0,$\emph{ i.e.} $\varphi(s)\neq t$ and therefore $U\not\subset\varphi^{-1}(\{t\}).$
\CQFD 
\begin{Rem}\label{rem1} Note that in a topological space $S$, and for a continuous map $\varphi:S\rightarrow S$, the following conditions are equivalent~:
\begin{enumerate}[i)]
\item for every $t\in S,\ \varphi^{-1}(\{t\})$ is nowhere dense in $S$
\item for every nonempty open subset $U$ of $S,\ \varphi(U)$ is infinite.
\end{enumerate}
\end{Rem}
Indeed, if there exists a nonempty open subset~$U$ of~$S$ such that $U\subset\varphi^{-1}(\{t\})$ then $\varphi$ is constant on~$U$, so $ii)\Rightarrow i).$ Moreover if $\varphi(U)=\{s_1,\ldots,s_n\}$ for an open subset $U$ of $S,\ n\geq1,$ then 
\begin{align*}
\{s\in U\mid\ \varphi(s)=s_1\}&=\{s\in U\mid\varphi(s)\neq s_k,\ 2\leq k\leq n\}\\
&=\varphi^{-1}\big(S\backslash\{s_2,\ldots,s_n\}\big)\cap U.
\end{align*} 
The set $S\backslash\{s_2,\ldots,s_n\}$ is open in $S$, so $\{s\in U\mid\ \varphi(s)=s_1\}$ is a nonempty open subset of~$U$ (and of~$S$) although by $i),\ \{s\in S\mid\ \varphi(s)=s_1\}$ must have empty interior.\par
The previous remark, Theorem~\ref{faiblcomp} and Theorem~\ref{rec1} give the following~: 
\begin{Cor}
 Let $S$ be a compact Hausdorff space without isolated point. Then $\|uC_\varphi+T\|=\|u\|_\infty+\|T\|$ for every weakly compact operator $T:C(S)\rightarrow C(S)$ if and only if $\vert u\vert$ is constant on $S$ and the set $\varphi^{-1}(\{t\})$ is nowhere dense in $S,$ for every $t\in S.$ 
\end{Cor}
\textbf{Application : a negative answer to a question of Popov~\cite{P}}\\ 
Note that if $\varphi$ is onto and $\vert u\vert=1$ then $uC_\varphi$ is an isometry on $C(S).$ In~\cite{P}, Popov shows that every into isometry $J:L^1([0,1])\rightarrow L^1([0,1])$ is a Daugavet center. He raises the question whether this result is true when we substitute $L^1([0,1])$ with a Banach space $X$ having the Daugavet property. Actually, this is not true for $X=C(S)$. To see this, consider any composition operator whose symbol $\varphi$ is onto and constant on a nonempty open subset of $S.$ Then $C_\varphi$ is an isometry but there exists rank-$1$ operators on~$C(S)$ which does not satisfy equation~$(E_{1,\varphi}).$ \\
After our work was completed, an example was independently produced in~\cite{BK}. The authors considered a weighted composition operator $uC_\varphi:C[0,1]\rightarrow C[0,1]$ whose symbol $\varphi$ is constant on~$]1/2,1]$ and whose weight has not constant modulus on $[0,1].$ 

\subsection{Convex combinations of composition operators}
One can wonder if the set of Daugavet centers is a convex set. Actually it is easy to see that this is not true in full generality. Indeed, consider $u(x)=e^{2i\pi x}$ and $v(x)=e^{-2i\pi x},\ x\in[0,1].$ Then $u,v\in C[0,1],\ \vert u\vert=\vert v\vert=1$ so $uI$ and $vI$ are Daugavet centers in~$C([0,1]),$ but $\big(u(x)+v(x)\big)/2=\cos 2\pi x$ which has not constant modulus on~$[0,1].$ Therefore $(uI+vI)/2$ is not a Daugavet center in $C[0,1].$ Nevertheless it turns out that a convex combination of particular (non zero) Daugavet centers can be a Daugavet center. Let us consider the case of composition operators.\\Note that a convex combination of composition operators is not in general a composition operator itself. Indeed, assume that $C_\varphi=tC_{\psi_1}+(1-t)C_{\psi_2}$ where $\varphi,\psi_1,\psi_2$ are continuous functions from~$S$ to~$S$ and $0<t<1.$ Assume that $\varphi\neq\psi_1$ and take $s_0$ such that $\varphi(s_0)\neq\psi_1(s_0).$ Now consider a open subset~$U$ of~$S$ such that $\varphi(s_0)\in U$ and $\psi_1(s_0)\notin U.$ Choose $f\in C(S),\ \|f\|_\infty=1$ satisfying $f\big(\varphi(s_0)\big)=1$ and $\vert f\vert<1$ out of~$U.$ Then\[1=\vert f\big(\varphi(s_0)\big)\vert\leq t\big\vert f\big(\psi_1(s_0)\big)\big\vert+(1-t)\big\vert f\big(\psi_2(s_0)\big)\big\vert<1\]
which leads to a contradiction.\par
\vspace{\baselineskip}
Let $\varphi$ and $\psi$ be continuous functions from $S$ to $S$. Assume that $\varphi\neq\psi.$ Define~$$S_1=\{s\in S\mid\ \varphi(s)\neq\psi(s)\}.$$ Then $S_1$ is a nonempty open subset of $S$ since $S$ has no isolated point. Consider convex combinations of $C_\varphi$ and $C_\psi.$ For $t\in[0,1],$we define $T_t=tC_\varphi+(1-t)C_\psi.$ Point out that~$\|T_t\|=1.$ For convenience, we note \[\Delta_T(s)=\big\vert t+\mu_s\big(\{\varphi(s)\}\big)\big\vert+\big\vert 1-t+\mu_s\big(\{\psi(s)\}\big)\big\vert-\Big(1+\big\vert\mu_s\big(\{\varphi(s)\}\big)\big\vert+\big\vert\mu_s\big(\{\psi(s)\}\big)\big\vert\Big),\] and \[\tilde{\Delta}_T(s)=\big\vert 1+\mu_s\big(\{\varphi(s)\}\big)\big\vert-\Big(1+\big\vert\mu_s\big(\{\varphi(s)\}\big)\big\vert\Big)\]
where~$(\mu_s)_{s\in S}$ is the family of measures representing~$T$ and $s\in S.$
As for weighted composition operators, we have the following property~:
\begin{Prop}\label{eq2conv}
Let $T$ be an operator on $C(S).$ Assume that the family of measures $(\mu_s)_{s\in S}$ representing~$T$ satisfies the condition~: 
for every nonempty open set$\ U\subset S$~:\par
-If $U\cap S_1\neq\emptyset,$ then
\begin{equation}\label{eq4}
\sup_{s\in U\cap S_1}\Delta_T(s)=0\end{equation}\par
-If $U\cap S_1=\emptyset,$ then
\begin{equation}\label{eq5}
\sup_{s\in U}\tilde{\Delta}_T(s)=0.
\end{equation}
Then the following equation holds true :  \[\|T_t+T\|=1+\|T\|.\]
\end{Prop}
{\bf Proof. }
One only has to consider open subsets $U$ of $S$ of the form $U=\{s\in S\mid\ \|\mu_s\|>\|T\|-\eps\}$ where $\eps>0.$ If $U\cap S_1=\emptyset$ then $\varphi=\psi$ on $U$ so the proof of Lemma~\ref{eq1} tells us that $\|T_t+T\|\geq\sup_{s\in U}\|\delta_{\varphi(s)}+\mu_s\|\geq 1+\|T\|-\eps.$ Else, $\|T_t+T\|\geq\sup_{s\in U\cap S_1}\|t\delta_{\varphi(s)}+(1-t)\delta_{\psi(s)}+\mu_s\|$ which is greater than $1+\|T\|-\eps$ using the same method as in the proof of Lemma~\ref{eq1}. \CQFD
From this we can deduce that any convex combination of composition operators which are Daugavet centers is still a Daugavet center.
\begin{Theo}\label{convex}
Assume that $C_\varphi$ and $C_\psi$ are Daugavet center. Then every weakly compact operator $T$ on~$C(S)$ satisfies the norm equation \[\|tC_\varphi+(1-t)C_\psi+T\|=1+\|T\|,\] for all $t\in[0,1].$
\end{Theo}
{\bf Proof. }Take $t\in[0,1].$ Argue by contradiction and assume that the family $(\mu_s)_{s\in S}$ does not satisfy conditions of Proposition~\ref{eq2conv}.\\
\emph{First case}~: Let $U$ be a nonempty open subset of~$S$ such that $U\cap S_1\neq\emptyset$ and \eqref{eq4} does not hold. Then there exists $\beta>0$ such that \[\big\vert t+\mu_s\big(\{\varphi(s)\}\big)\big\vert+\big\vert 1-t+\mu_s\big(\{\psi(s)\}\big)\big\vert-\Big(1+\big\vert\mu_s\big(\{\varphi(s)\}\big)\big\vert+\big\vert\mu_s\big(\{\psi(s)\}\big)\big\vert\Big)<-4\beta\]
for every $s\in U\cap S_1.$ Then \[\big\vert\mu_s\big(\{\varphi(s)\}\big)\big\vert+\big\vert\mu_s\big(\{\psi(s)\}\big)\big\vert>2\beta,\quad\forall s\in U\cap S_1.\]
Let \begin{align*}
V_1&=\{s\in U\cap S_1\mid\ \big\vert\mu_s\big(\{\varphi(s)\}\big)\big\vert>\beta\}\\
V_2&=\{s\in U\cap S_1\mid\ \big\vert\mu_s\big(\{\psi(s)\}\big)\big\vert>\beta\}.
\end{align*}
Since $U\cap S_1\subset V_1\cup V_2$, we can assume without loss of generality that $V_1$ contains a nonempty open set $V.$ So $\big\vert\mu_s\big(\{\varphi(s)\}\big)\big\vert>\beta$ for every $s\in V.$ Then we follow the proof of Theorem~\ref{faiblcomp} to obtain a contradiction.\\
\emph{Second case}~: If $U$ is a nonempty open subset of~$S$ such that $U\subset S\backslash S_1$ and \eqref{eq5} does not hold, then the same proof as in Theorem~\ref{faiblcomp} leads to a contradiction.\CQFD   

\subsection{Operators factorizing through an Asplund space}
The aim of this section is to extend a result of Ansari in~\cite{A} stating that every operator on $C(S)$ factorizing through $c_0$ satisfies the Daugavet equation. 
Let $T:C(S)\rightarrow C(S)$ be an operator, where $S$ is a compact Hausdorff space, and $(\mu_s)_{s\in S}$  the family of measures associated to~$T.$ Note that if $\mu_s\big(\{\varphi(s)\}\big)=0$ for all $s\in S$, and $\vert u\vert$ is constant, then $T$ trivially satisfies condition $(2.2)$ of corollary~(\ref{eq2}). Actually, it is sufficient that the measures $(\mu_s)$ almost satisfies this condition~: \par Define $S_\eps=\big\{s\in S\mid\ \big\vert\mu_s\big(\{\varphi(s)\}\big)\big\vert<\eps\big\},$ for each $\eps>0.$ We have the following~:
\begin{Lem}\label{dens}
If $\vert u\vert$ is constant and if the sets $S_\eps$ are dense in $S,$ for every $\eps>0,$ then \[\|uC_\varphi+T\|=\|u\|_\infty+\|T\|.\]
\end{Lem}
{\bf Proof.} Take $U$ a nonempty open set in $S$ and $\eps>0.$ By density of $S_\eps$ in $S$, there exists  $s_\eps\in U$  satisfying $\big\vert\mu_{s_\eps}\big(\{\varphi(s_{\eps})\}\big)\big\vert<\eps,$ and so
\begin{align*}
\big\vert u(s_\eps)+\mu_{s_\eps}\big(\{\varphi(s_\eps)\}\big)\big\vert-\big(\|u\|_\infty+\big\vert\mu_{s_\eps}(\{\varphi(s_\eps)\})\big\vert\big) & \geq-2\big\vert\mu_{s_\eps}(\{\varphi(s_\eps)\})\big\vert \\
 &>-2\eps.
\end{align*}
Thus for every nonempty open set $U$ of $S$, we have\[\dis\sup_{s\in U}\bigg(\big\vert u(s)+\mu_s\big(\{\varphi(s)\}\big)\big\vert-\Big(\|u\|_\infty+\big\vert\mu_s\big(\{\varphi(s)\}\big)\big\vert\Big)\bigg)=0.\]
We conclude with corollary ~\eqref{eq2}.\CQFD
Now we can prove the following result : 
\begin{Theo}\label{dualsep}
Let $S$ be a compact Hausdorff space without isolated point, $\varphi:S\rightarrow S$ a continuous map,~$u\in C(S)$ and $T:C(S)\rightarrow C(S)$ an operator such that $T^*(M(S))$ is separable. If $\varphi^{-1}(\{t\})$ is nowhere dense in $S$, for every $t\in S,$ and if $\vert u\vert$ is constant, then $T$ satisfies equation $\|uC_\varphi+T\|=\|u\|_\infty+\|T\|.$
\end{Theo}
{\bf Proof. }Let $\{\rho_n,\ n\in\N\}$ be a dense subset of $T^*\big(M(S)\big).$ As previously, $S_\eps=\big\{s\in S\mid\ \big\vert\mu_s\big(\{\varphi(s)\}\big)\big\vert<\eps\big\},$ where $\mu_s=T^*(\delta_s),$ and $A=\dis\bigcap_{n\geq0}\big\{s\in S\mid\ \rho_n\big(\{\varphi(s)\}\big)=0\big\}.$ We want to show that~:
\begin{enumerate}[i)]
\item $A$ is dense in $S.$
\item $\forall\eps>0,\ A\subset S_\eps.$
\end{enumerate}
Then we conclude with Lemma~\ref{dens}.\\
To prove $i),$ we are going to show that $S\backslash A$ is nowhere dense. Indeed, \[S\backslash A=\bigcup_{n\geq0}\big\{s\in S\mid\ \rho_n\big(\{\varphi(s)\}\big)\neq0\big\}=\bigcup_{n\geq0}\bigcup_{p\geq1}A_{n,p}\]
where $A_{n,p}=\big\{s\in S\mid\ \big\vert\rho_n\big(\{\varphi(s)\}\big)\big\vert>\frac{1}{p}\big\}.$ Since $\rho_n$ is a finite measure, this implies that the sets $\varphi(A_{n,p})$ are finite (hence closed) for every $n\geq0$, $p\geq1.$ But $A_{n,p}\subset\varphi^{-1}\big(\varphi(A_{n,p})\big)$ which a finite union of nowhere dense sets (c.f. Remark~\ref{rem1}). Using Baire's theorem, $S\backslash A$ is contained in a nowhere dense set, and $A$ is dense in $S$.\\
Proof of $ii)$ : let $s\in S$ and $\eps>0.$ By density of $(\rho_n)_n$ in $T^*(M(S))$, there exists an integer $n_0\geq0$ such that $\|T^*(\delta_s)-\rho_{n_0}\|<\eps.$ Then $\vert\mu_s-\rho_{n_0}\vert\big(\{\varphi(s)\}\big)<\eps.$ Taking $s\in A,$ it follows that $\big\vert\mu_s\big(\{\varphi(s)\}\big)\big\vert<\eps,$ \textit{i.e.} $A\subset S_\eps.$
\CQFD 
If $T:C(S)\rightarrow C(S)$ factorizes through a space~$X$ having a separable dual, then Theorem~\ref{dualsep} applies. In particular this holds for the class of operators factorizing through $c_0.$ Actually, regarding operators factorizing through a space~$X$, one does not need to assume that $X^*$ is separable in the case where $S$ is metrizable. We recall the following definition~:
\begin{Def}
A Banach space $X$ is called an Asplund space if its dual space has the Radon-Nikod\'ym property.
\end{Def}
Every dual space which is separable has the Radon-Nikod\'ym property, and so every Banach space with separable dual is Asplund. Asplund spaces are characterized by the fact that every separable subspace has a separable dual.\par
\begin{Cor}\label{asplund}
Let $S$ be a metric compact space without isolated points, $\varphi:S\rightarrow S$ a continuous map,~$u\in C(S)$ and $T:C(S)\rightarrow C(S)$ an operator factorizing through an Asplund space~$X$.  If $\varphi^{-1}(\{t\})$ is nowhere dense in $S$, for every $t\in S,$ and if $\vert u\vert$ is constant on~$S$, then $T$ satisfies equation $\|uC_\varphi+T\|=\|u\|_\infty+\|T\|.$
\end{Cor}
{\bf Proof. }Write $T=T_2T_1$ with $T_1:C(S)\rightarrow X$ and $T_2:X\rightarrow C(S).$ Since $S$ is a  metric compact space, $C(S)$ is separable. So we can assume, by replacing $X$ by $\overline{T_1\big(C(S)\big)}$ that $X^*$ is separable. Thus $T^*(M(S))$ is separable, and the result follows from Theorem~\ref{dualsep}. \CQFD
\begin{Rem}\label{rem2} Since every compact operator factorizes through a subspace of $c_0$, this gives another proof of Theorem~\ref{faiblcomp} for compact operators on~$C(S).$ Moreover every weakly compact operator factorizes through a reflexive space (which is Asplund), giving another proof of Theorem~\ref{faiblcomp} for weakly compact operators on $C(S)$ where $S$ is a metric compact space without isolated point.
\end{Rem}
In the case where $uC_\varphi=I$, Theorem~\ref{dualsep} is a particular case of an already known result in Banach spaces with the Daugavet property. If we consider a Banach space~$X$ having the Daugavet property, then every operator~$T:X\rightarrow X$ such that $T^*(X^*)$ is separable satisfies the Daugavet equation. This can be seen by using a result of Shvidkoy~\cite{Sh} which says that an operator $T:X\rightarrow X$ not fixing a copy of~$\ell_1$ satisfies the Daugavet equation. Then it is obvious that if~$T$ fixes a copy of~$\ell_1$ then~$T^*$ fixes a copy of~$\ell_\infty$, hence~$T^*(X^*)$ is not separable.\par
As another immediate consequence of Theorem~\ref{dualsep}, we have the following for particular weighted composition operators (which can also be viewed directly)~:
\begin{Cor}
Let $S$ be a compact Hausdorff space without isolated points, $\varphi:S\rightarrow S$ be a continuous map and $u\in C(S).$
If for every $t\in S,\ \varphi^{-1}(\{t\})$ is nowhere dense in $S,$ and if $\vert u\vert$ is constant, then $uC_\varphi:C(S)\rightarrow C(S)$ does not factorize through a space having a separable dual space. If~$S$ is metrizable, then $uC_\varphi$ does not factorize through an Asplund space.
\end{Cor}

\section{Equation $(E_{u,\varphi})$ for classes of operators on $A(\D)$}

In this section, we want to adapt D. Werner's method in~\cite{We3} to find new Daugavet centers in subspaces of $C(S)$-spaces, and particularly for the disk algebra $A(\D).$ Actually we will consider weighted composition operators $uC_\varphi$ on a functionnal Banach space~$X$ and will formulate conditions on an isometric embedding of~$X$ into~$C(S)$ implying that~$X$ is $(u,\varphi)$-nicely embedded. Then we find conditions so that every weakly compact operator on a $(u,\varphi)$-nicely embedded space satisfies equation $\|uC_\varphi+T\|=\|uC_\varphi\|+\|T\|.$
\subsection{General approach}
Let $(X,\|.\|)$ denotes a functional Banach space on~$\Omega\ \big(X\subset\mathcal{F}(\Omega,\C)\big).$ Consider $\varphi$ a map such that $\varphi(\Omega)\subset\Omega$ and $u\in X$ such that $0<\|u\|<\infty.$ Assume that $uC_\varphi:f\in X\mapsto u.(f\circ\varphi)\in X$ is a weighted composition operator acting continuously on~$X.$
Let $S$ be a compact Hausdorff space without isolated point. An isometry $J:X\rightarrow C(S)$ is said to be a~$(u,\varphi)$-\textit{nice embedding} and ~$X$ is said to be~$(u,\varphi)$-\textit{nicely embedded} into $C(S)$ if the following conditions are satisfied for every $s\in S$~:
\begin{description}
\item[(C1)] if $p_s=(uC_\varphi)^*J^*(\delta_s)\in X^*$, then $\|p_s\|=\|u\|>0.$
\item[(C2)] $\mathrm{Vect}(p_s)$ is an $L$-summand in $X^*.$
\end{description}
Recall that a closed subspace $F$ of a Banach space $E$ is an $L$-summand if there exists a projection $\Pi$ from $E$ onto $F$ such that, for every $x\in E,$  \[\|x\|=\|\Pi x\|+\|x-\Pi x\|.\]
We say that $F$ is an $M$-ideal if its annihilator $F^\bot\subset E^*$ is an $L$-summand. Then condition \textbf{(C2)} can be reformulated as~: $\ker(p_s)$ is an $M$-ideal in~$X$. Condition \textbf{(C1)} forces $\|uC_\varphi\|=\|u\|.$\par
Assume that $X$ is $(u,\varphi)$-nicely embedded in $C(S).$ Condition \textbf{(C2)} provides us a family of projections $(\Pi_s)_{s\in S}$ satisfying\[\|x^*\|=\|\Pi_sx^*\|+\|x^*-\Pi_sx^*\|,\quad \textrm{for every } x^*\in X^*\] and a family $(\pi_s)_{s\in S}$ in $X^{**}$ such that \[\Pi_sx^*=\pi_s(x^*)p_s,\quad\textrm{for every } x^*\in X^*.\]
Note that $\pi_s(p_s)=1.$\par
Consider the equivalence relation~$\sim$ on $S$\[s\sim t\Leftrightarrow \Pi_s=\Pi_t.\] Note $E_s$ the class of~$s$ in $S$. Then $E_s$ is closed, and condition \textbf{(C1)} tells us that $E_s=\{t\in S\mid\ p_t=\lambda p_s,\ \lambda\in\T\}.$ We will need the following condition~:
\begin{description}
\item[(C3)] for all $s\in S,$ the class $E_s$ is nowhere dense in $S.$
\end{description}
Let $T:X\rightarrow X$ be an operator, and $q_s=(JT)^*(\delta_s)\in X^*,\ s\in S.$ Then $s\mapsto q_s$ is continuous for the weak$^*$-topology on $X^*,$ and $\|T\|=\sup_S\|q_s\|.$\\
We can now express some results, whose proofs are similar to those in section~$2$ and are given in ~\cite{We3} in the particular case where $\varphi(x)=x,\ x\in\Omega$ and $u\equiv1.$
\begin{Prop}\label{prop1}
Suppose $X$ is $(u,\varphi)$-nicely embedded in $C(S),$ and $T$ is an operator acting on~$X.$
Then \[\|uC_\varphi+T\|=\|u\|+\|T\|\] if and only if  
\[for\ every\ \eps>0,\quad \sup_{\{s \mid\ \|q_s\|>\|T\|-\eps\}}\Big(\vert1+\pi_s(q_s)\vert-\big(1+\vert\pi_s(q_s)\vert\big)\Big)=0.\]
\end{Prop}
\begin{Prop}\label{eq3}
Suppose that $X$ is $(u,\varphi)$-nicely embedded in $C(S)$, and that condition \textbf{(C3)} holds. Let $T$ be an operator on $X.$ If we have
\[\textrm{for all }t\in S,\ s\mapsto\pi_t(q_s)\textrm{ is continuous,}\]
then $T$ satisfies equation $(E_{u,\varphi}).$
\end{Prop}
\begin{Rem}Every weakly compact operator $T$ on~$X$ fulfills conditions  of Proposition~\ref{eq3}, and consequently the equality $\|uC_\varphi+T\|=\|u\|_\infty+\|T\|$ holds.
\end{Rem}
We want to obtain this result for the class of operators whose adjoint has separable range. Let us start with a lemma which will be useful for the proof of the next proposition~:
\begin{Lem}\label{lemnoneq}
(\cite{We3}, Lemma $2.3$)\quad Suppose $X$ is $(u,\varphi)$-nicely embedded in $C(S).$ If $t_1,\dots,t_k$ are pairwise nonequivalent points (for the relation $\sim$), then\[\|x^*\|\geq\sum_{j=1}^k\|\Pi_{t_j}(x^*)\|,\quad\textrm{for every } x^*\in X^*.\]
\end{Lem}
\begin{Prop}\label{dualsep2}
Let $X$ be a $(u,\varphi)$-nicely embedded space in $C(S)$ and satisfying condition \textbf{(C3)}, and $T$ be an operator on~$X$ such that $T^*(X^*)$ is separable.
Then \[\|uC_\varphi+T\|=\|u\|+\|T\|\]
\end{Prop}
{\bf Proof. }Consider the sets $S_\eps=\{s\in S\mid\ \vert\pi_s(q_s)\vert<\eps\},$ where $q_s=(JT)^*(\delta_s)$ and $\eps>0.$
If we show that $S_\eps$ is dense in $S,$ for every $\eps>0,$ then $T$ fulfills the condition of Proposition~\ref{prop1}.\par
Let $\{\psi_n,\ n\in\N\}$ be a dense subset of $T^*(X^*)$ and define $$A=\big\{s\in S\mid\ \pi_s(\psi_n)=0,\ \forall n\in \N\big\}.$$ As in the proof of Theorem~\ref{dualsep} we want to show that~:
\begin{enumerate}[i)]
\item $A$ is dense in $S$
\item $\forall\eps>0,\ A\subset S_\eps.$
\end{enumerate}
The proof of $ii)$ is similar to the one in Theorem~\ref{dualsep}. For $i),$ we show that $S\backslash A$ is nowhere dense. Indeed \[S\backslash A=\bigcup_{n\geq0}\big\{s\in S\mid\ \pi_s(\psi_n)\neq0\big\}=\bigcup_{n\geq0}\bigcup_{p\geq1}A_{n,p}\]
where $A_{n,p}=\big\{s\in S\mid\ \big\vert\pi_s(\psi_n)\big\vert>\frac{1}{p}\big\}.$ By Lemma~\ref{lemnoneq} there is a finite number of equivalence classes for $\sim$ in $A_{n,p}$ (less than $p\|\psi_n\|/\|u\|$). These equivalence classes are closed and nowhere dense (by condition \textbf{C3}). The Baire property yields that $S\backslash A$ is contained in a nowhere dense set, implying that $A$ is dense in~$S.$
\CQFD 
In the case where~$X$ is separable, we have the same result as in Corollary~\ref{asplund}.
\begin{Cor}\label{asplund2}
Let $X$ be a $(u,\varphi)$-nicely embedded space in $C(S)$ satisfying condition \textbf{(C3)}, and $T$ be an operator on~$X$ which factorizes through an Asplund space~$E$. Assume that~$X$ is separable.
Then \[\|uC_\varphi+T\|=\|u\|+\|T\|\]
\end{Cor}
\subsection{Applications}
Obviously one can apply this results to the case $X=C(S)$ where~$S$ is a compact \mbox{Hausdorff} space without isolated point, with~$J$ the natural inclusion into~$C(S).$ If~\mbox{$\varphi:S\rightarrow S$} is continuous and if~$u\in C(S)$ with $u\neq0,$ then~$p_s=u(s)\delta_{\varphi(s)}$ so that condition \textbf{(C1)} forces~$\vert u\vert$ to be constant equal to $\|u\|_\infty.$ Then \mbox{$\Pi_s:\mu\in C(S)^*\mapsto\big(\mu\big(\{\varphi(s)\}\big)/u(s)\big)p_s$} is a $L$-projection, and condition~\textbf{(C2)} holds. Finally for $s$ and $t$ in $S,$ 
\begin{align*}
s\sim t &\Leftrightarrow\exists\lambda\in\T,\ p_s=\lambda p_t \\
                     &\Leftrightarrow\exists\lambda\in\T,\ \vert u(s)\vert\delta_{\varphi(s)}=\lambda\vert u(t)\vert\delta_{\varphi(t)}\\
                     &\Leftrightarrow\varphi(s)=\varphi(t).
\end{align*}
Thus $E_s=\varphi^{-1}\big(\{\varphi(s)\}\big).$ So condition \textbf{(C3)} is fulfilled if $\varphi^{-1}(\{t\})$ is nowhere dense, for every $t\in S.$ Therefore we recover most of the results of section~$2.$
\vspace{\baselineskip}\par
We now turn to the disk algebra $A(\D).$  Let $\D=\{z\in \C\mid\ \vert z\vert<1\}$ denote the unit disk. The disk algebra $A(\D)$ is the algebra of holomorphic maps on $\D$ which are continuous on $\overline{\D}$, endowed with the supremum norm $\|f\|_\infty=\sup\{\vert f(z)\vert\mid\ z\in\D\}.$ Considering $u\neq0$ and $\varphi$ in the disk algebra with $\|\varphi\|_\infty\leq1,$ one can define the weighted composition operator $uC_\varphi$ acting on $A(\D),$ with $\|uC_\varphi\|=\|u\|_\infty.$ We will assume that $\varphi$ is not constant (implying $\varphi(\D)\subset\D$), otherwise $uC_\varphi$ is a rank-$1$ operator and therefore is not a Daugavet center. On the other hand, $C_\varphi$ is compact on~$A(\D)$ if and only if $\|\varphi\|_\infty<1.$ Hence a necessary condition ensuring that every weakly compact operator on $A(\D)$ fulfills equation $(E_{u,\varphi})$ is that $\|\varphi\|_\infty=1.$ Actually, we are going to prove that a strong negation of ``$\|\varphi\|_\infty<1$'' is necessary.\par
Consider the isometry \mbox{$J:f\in A(\D)\mapsto J(f)=f_{|_\T}\in C(\T).$} It is well known that the image of $A(\D)$ by $J$ is the closed space $\{f\in C(\T)\mid\ \hat{f}(n)=0,\ \forall n<0\}.$ We want conditions on $\varphi$ and $u$ implying that $A(\D)$ is $(u,\varphi)$-nicely embedded in $C(\T)$, and moreover that condition \textbf{(C3)} is fulfilled.\par
For $\omega\in\T,$ let\[\ p_\omega:=(uC_\varphi)^*J^*(\delta_\omega)=u(\omega)\delta_{\varphi(\omega)_{|_{A(\D)}}}\in A(\D)^*.\]Clearly $\|p_\omega\|=\vert u(\omega)\vert,$ so \textbf{(C1)} is fulfilled if and only if $\vert u\vert$ is constant on $\T.$ Assume that $\vert u\vert$ is constant on $\T.$ To check condition \textbf{(C2)}, we have to show that $\ker(p_\omega)$ is an $M$-ideal. But 
\begin{align*}
\ker(p_\omega)&=\{f\in A(\D)\mid\ u(\omega)f\big(\varphi(\omega)\big)=0\}\\
                      &=\{f\in A(\D)\mid\ f\big(\varphi(\omega)\big)=0\}
\end{align*}
since $u(\omega)\neq 0.$ It is an $M$-ideal if and only if $\varphi(\omega)\in\T$ (see~\cite{HWW} p.~4). It means that \textbf{(C2)} is fulfilled if $\varphi$ is an inner function. Finally, if $\omega_1,\omega_2\in\T,$ 
\begin{align*}
\omega_1\sim\omega_2 &\Leftrightarrow\exists\lambda\in\T,\ p_{\omega_1}=\lambda p_{\omega_2} \\
                     &\Leftrightarrow\exists\lambda\in\T,\ \vert u(\omega_1)\vert\delta_{\varphi(\omega_1)}=\lambda\vert u(\omega_2)\vert\delta_{\varphi(\omega_2)}\textrm{ on }A(\D)\\
                     &\Leftrightarrow\varphi(\omega_1)=\varphi(\omega_2).
\end{align*}
Thus $E_\omega=\varphi^{-1}\big(\{\varphi(\omega)\}\big)\cap\T.$ If $\varphi$ is not constant, then condition \textbf{(C3)} is fulfilled.\par
To summarize, we have the following~:
\begin{Prop}
Let $\varphi$ be an inner function and $u$ be a multiple of an inner function. If $T:A(\D)\rightarrow A(\D)$ is such that $T^*\big(A(\D)^*\big)$ is separable, then equation $\|uC_\varphi+T\|=\|u\|_\infty+\|T\|$ holds true.
\end{Prop}
Since $A(\D)$ is separable, Remark~\ref{rem2} and Corollary~\ref{asplund2} gives~:
\begin{Cor}\label{prop}
Let $\varphi$ be an inner function and $u$ be a multiple of an inner function. Then every weakly compact operator~$T:A(\D)\rightarrow A(\D)$ satifies equation$$\|uC_\varphi+T\|=\|u\|_\infty+\|T\|.$$
\end{Cor}
This corollary leads to the following remark on (general) essential norms of weighted composition operators on the disk algebra.
\begin{Rem}
Let $X$ be a Banach space, $B(X)$ be the space of bounded operator on~$X,$ $\mathcal{K}(X)$ be the closed subspace of $B(X)$ consisting of compact operators on~$X$ and $\mathcal{W}(X)$ be the closed subspace of $B(X)$ consisting of weakly-compact operators on~$X.$ Recall that if $\mathcal{I}$ is a closed subspace of~$B(X),$ the essential norm (relatively to~$\mathcal{I}$) of $S\in B(X)$ is the distance from $S$ to $\mathcal{I}$~:\[\|S\|_{e,\mathcal{I}}=\inf\{\|S+T\|;\ T\in \mathcal{I}\}.\]
This is the canonical norm on the quotient space~$B(X)/\mathcal{I}$. The classical case corresponds to the case of compact operators $\mathcal{I}=\mathcal{K}(X)$. In this case, the above quotient space is the Calkin algebra. General essential norms of weighted composition operators on $A(\D)$ are estimated in~\cite{Le}. When $\mathcal{I}\subset\mathcal{W}\big(A(\D)\big),$ and in the particular case where $\varphi$ is an inner function $(\varphi(\D)\subset \D)$ and $u$ is a multiple of an inner function, Corollary~\ref{prop} not only gives us the essential norm relatively to $\mathcal{I}$ of~$uC_\varphi,$ but how the norm of $uC_\varphi$ reacts under perturbation by operators in the class~$\mathcal{I}.$\par
\end{Rem}
Although D. Werner's method gives sufficient conditions ensuring that weighted composition operators are Daugavet centers, it turns out that these ones are also necessary. 
\begin{Prop}
Suppose that every rank-$1$ operator on $A(\D)$ satisfies equation $(E_{u,\varphi}).$ Then $\varphi$ is inner and  $\vert u\vert$ is constant on $\T.$
\end{Prop}
{\bf Proof. }Assume $\varphi$ is not inner. Then there exists $\omega\in\T$ such that $\vert\varphi(\omega)\vert=r<1.$ As $u$ is not constant equal to zero, we can assume, taking if necessary a $\omega'\in\T$ close to $\omega,$ that $u(\omega)\neq0.$ Let $g\in A(\D)$ defined by $g(z)=(1+\bar{\omega}z)/2,$ where $z\in\D.$ We consider the rank-$1$ operator $T:f\mapsto Tf=u(\omega)f\big(\varphi(\omega)\big)g,$ for all $f\in A(\D).$ We have $\|T\|=\vert u(\omega)\vert.$ For \mbox{$0<\eps<\min(1-r,\vert u(\omega)\vert/3),$} there exists an arc $I_\omega\subset\T$ containing $\omega$ such that for every $z\in I_\omega$, we have $\vert \varphi(z)-\varphi(\omega)\vert\leq\eps,\ \vert1-g(z)\vert<1/2$ and $\vert u(z)-u(\omega)\vert<\eps.$ Let $f\in A(\D)$ with $\|f\|_\infty=1$~:   
\[\|uC_\varphi(f)-Tf\|=\sup_{\vert z\vert=1}\big\vert u(z)f\big(\varphi(z)\big)-u(\omega)f\big(\varphi(\omega)\big)g(z)\big\vert\]
If $z\notin I_\omega$, then $\big\vert u(z)f\big(\varphi(z)\big)-u(\omega)f\big(\varphi(\omega)\big)g(z)\big\vert\leq                 
\|u\|_\infty+\vert u(\omega)\vert\dis\sup_{z\in\T\backslash I_\omega}\vert g(z)\vert.$\\ For any $a,b\in D(0,r+\eps)$, we have by the Cauchy formula :
\begin{align*}
\vert f(a)-f(b)\vert &=\bigg\vert\frac{1}{2i\pi}\int_{\vert z\vert=1}\frac{a-b}{(z-a)(z-b)}f(z)\ud z\bigg\vert \\
                     &\leq\frac{\vert a-b\vert}{2\pi}\int_0^{2\pi}\frac{\big\vert f\big(e^{i\theta}\big)\big\vert}{\big\vert e^{i\theta}-a\big\vert\big\vert e^{i\theta}-b\big\vert}\ud\theta \\
                     &\leq\frac{\vert a-b\vert}{\big(1-(r+\eps)\big)^2}.
\end{align*}
Thus for $z\in I_\omega,\ \varphi(z)\in D(0,r+\eps)$ and we have~: 
\begin{align*}
\big\vert u(z)f\big(\varphi(z)\big)-u(\omega)f\big(\varphi(\omega)\big)g(z)\big\vert  & \leq \vert u(z)\vert\big\vert f\big(\varphi(z)\big)-f\big(\varphi(\omega)\big)\big\vert+\vert u(z)-u(\omega)\vert\big\vert f\big(\varphi(\omega)\big)\big\vert\\
& \quad+\big\vert u(\omega)f\big(\varphi(\omega)\big)\big\vert\vert1-g(z)\vert\\
  &\leq  \|u\|_\infty\bigg(\frac{\eps}{\big(1-(r+\eps)\big)^2}\bigg)+\eps+\frac{\vert u(\omega)\vert}{2}\\
  &\leq  \|u\|_\infty+\frac{5}{6}\vert u(\omega)\vert
\end{align*}
for a suitable $\eps>0.$ So
\[\|uC_\varphi(f)-Tf\|\leq\max\bigg(\|u\|_\infty+\frac{5}{6}\vert u(\omega)\vert,\|u\|_\infty+\vert u(\omega)\vert\delta\bigg),\]
where $\delta=\dis\sup_{z\in\T\backslash I_\omega}\vert g(z)\vert<1.$ This gives $\|C_\varphi-T\|<\|u\|_\infty+\vert u(\omega)\vert=\|C_\varphi\|+\|T\|$ which is absurd. So $\varphi$ is an inner function.\par
We use a similar argument to show that $\vert u\vert$ is constant on the unit circle.
\CQFD
We summarize our results in the following corollary~: 
\begin{Cor}
Let $\varphi\in A(\D),\ \varphi(\D)\subset\D$ and $u\in A(\D).$ Then $uC_\varphi$ is a Daugavet center in $A(\D)$ if and only if $\varphi$ is an inner function and $u$ is a multiple of an inner function.
\end{Cor}
\begin{Rem}
The case of the disk algebra is different from the case of $C(S)$. Indeed, we have seen that a function $\varphi:S\rightarrow S$ could induce a composition operator $C_\varphi$ on $C(S)$ which is an isometry but is not a Daugavet center (see ~section~$2$). Whereas if $\varphi\in A(\D)$ satisfies $\varphi(\D)\subset\D$, $\varphi$ is an inner function if and only if $C_\varphi$ is an isometry, and so $C_\varphi$ is an isometry if and only if $C_\varphi$ is a Daugavet center. 
\end{Rem}
Note that in the particular case where $\varphi$ is a disk automorphism it is easy to see that $C_\varphi$ is a Daugavet center. Indeed, consider a weakly compact operator~$T$ on~$A(\D).$ Then $C_\varphi+T=C_\varphi(I+C_\varphi^{-1}T).$ The fact that $C_\varphi$ is an isometry implies that $\|C_\varphi+T\|=\|I+C_\varphi^{-1}T\|.$ Since the disk algebra has the Daugavet property (see~\cite{Wo}) and $C_\varphi^{-1}$ is itself an isometry, we have \[\|C_\varphi+T\|=1+\|C_\varphi^{-1}T\|=1+\|T\|.\]

\vspace{\baselineskip}
Romain Demazeux, \it{Univ Lille Nord de France F-59 000 LILLE, FRANCE\\
UArtois, Laboratoire de Math\'ematiques de Lens EA 2462,\\
F\'ed\'eration CNRS Nord-Pas-de-Calais FR 2956,\\
F-62 300 LENS, FRANCE}\\
romain.demazeux@euler.univ-artois.fr

\end{document}